%% This is an Amstex file, orginal file `tanjanov.tex'. It is the version
%% of October 26, 2005, not for publication
%% or circulation!!Minor changes December 13, 2005,
%% and December 25, 29,  2005, April 16, 2006.
%% Corrections 30.10.2010
%% Published "Publs. Inst. Math." 80(94) (2006), 141-156.
\input amstex.tex
\documentstyle{amsppt}
\magnification1200
 \hsize=12.5cm
 \vsize=18cm
\hoffset=1cm
\voffset=2cm

\def\DJ{\leavevmode\setbox0=\hbox{D}\kern0pt\rlap
{\kern.04em\raise.188\ht0\hbox{-}}D}
\def\dj{\leavevmode
 \setbox0=\hbox{d}\kern0pt\rlap{\kern.215em\raise.46\ht0\hbox{-}}d}

\def\txt#1{{\textstyle{#1}}}
\baselineskip=13pt
\def\hf{{\textstyle{1\over2}}}
\def\a{\alpha}\def\b{\beta}
\def\d{{\,\roman d}}
\def\e{\varepsilon}

\def\G{\Gamma}

\def\s{\sigma}
\def\t{\theta}
\def\={\;=\;}
\def\zx{\zeta(\hf+ix)}
\def\zt{\zeta(\hf+it)}

\def\D{\Delta}

\def\R{\Re{\roman e}\,} 
\def\z{\zeta}

 \def\t{\theta}
\def\hf{{\textstyle{1\over2}}}
\def\txt#1{{\textstyle{#1}}}

\def\Z{{\Cal Z}}
%%%%%%%%%%% Fonts macros %%%%%%%%%%%%
\font\tenmsb=msbm10
\font\sevenmsb=msbm7
\font\fivemsb=msbm5
\newfam\msbfam
\textfont\msbfam=\tenmsb
\scriptfont\msbfam=\sevenmsb
\scriptscriptfont\msbfam=\fivemsb
\def\Bbb#1{{\fam\msbfam #1}}

\def \NN {\Bbb N}
\def \CC {\Bbb C}
\def \RR {\Bbb R}
\def \ZZ {\Bbb Z}

\font\ff=cmr8
\def\txt#1{{\textstyle{#1}}}
\baselineskip=13pt

\font\teneufm=eufm10
\font\seveneufm=eufm7
\font\fiveeufm=eufm5
\newfam\eufmfam
\textfont\eufmfam=\teneufm
\scriptfont\eufmfam=\seveneufm
\scriptscriptfont\eufmfam=\fiveeufm
\def\mathfrak#1{{\fam\eufmfam\relax#1}}

\font\tenmsb=msbm10
\font\sevenmsb=msbm7
\font\fivemsb=msbm5
\newfam\msbfam
     \textfont\msbfam=\tenmsb
      \scriptfont\msbfam=\sevenmsb
      \scriptscriptfont\msbfam=\fivemsb
\def\Bbb#1{{\fam\msbfam #1}}

\def \NN {\Bbb N}
\def \CC {\Bbb C}

\def \RR {\Bbb R}
\def \ZZ {\Bbb Z}

  \def\rightheadline{{\hfil{\ff
Convolutions and mean square estimates
}\hfil\tenrm\folio}}

  \def\leftheadline{{\tenrm\folio\hfil{\ff
   Aleksandar Ivi\'c }\hfil}}
  \def\emptyheadline{\hfil}
  \headline{\ifnum\pageno=1 \emptyheadline\else
  \ifodd\pageno \rightheadline \else \leftheadline\fi\fi}

\font\ff=cmr8
\font\teneufm=eufm10
\font\seveneufm=eufm7
\font\fiveeufm=eufm5
\newfam\eufmfam
\textfont\eufmfam=\teneufm
\scriptfont\eufmfam=\seveneufm
\scriptscriptfont\eufmfam=\fiveeufm
\def\mathfrak#1{{\fam\eufmfam\relax#1}}

\font\tenmsb=msbm10
\font\sevenmsb=msbm7
\font\fivemsb=msbm5
\newfam\msbfam
\textfont\msbfam=\tenmsb
\scriptfont\msbfam=\sevenmsb
\scriptscriptfont\msbfam=\fivemsb
\def\Bbb#1{{\fam\msbfam #1}}

\def \NN {\Bbb N}
\def \CC {\Bbb C}

\def \RR {\Bbb R}
\def \ZZ {\Bbb Z}

\def\D{\Delta}
\def\a{\alpha}
\def\b{\beta} \def\e{\varepsilon}
 \def\d{\,{\roman d}}
%\topglue1cm
\topmatter
\title
Convolutions and mean square estimates of certain number-theoretic error terms
\endtitle
\author
Aleksandar Ivi\'c
\endauthor
\address
Katedra Matematike RGF-a, Universitet u Beogradu,  \DJ u\v sina 7,
11000 Beograd, Serbia and Montenegro
\endaddress
\subjclass 11 N 37, 11 M 06, 44 A 15, 26 A 12
\endsubjclass
\email {\tt  ivic\@rgf.bg.ac.rs, aivic\@matf.bg.ac.rs}
\endemail
\dedicatory
Dedicated to the memory of Tatjana Ostrogorski
\smallskip
``Publs. Inst. Math." 80(94) (2006), 141-156
\enddedicatory
\abstract
We study the convolution function
$$
C[f(x)] := \int_1^x f(y)f({x\over y}){\d y\over y}
$$
when $f(x)$ is a suitable number-theoretic error term.
Asymptotics and upper bounds for $C[f(x)]$ are derived
from  mean square bounds for $f(x)$. Some applications
are given, in particular to $|\zx|^{2k}$ and the classical Rankin--Selberg
problem from analytic number theory.
\endabstract
\endtopmatter
\document
\head
1. Convolution functions
\endhead
Motivated by considerations from analytic number theory,
the author investigated in [10] the following class of
convolution functions. Let ${\Cal M}_a$ denote
the set of functions $f(x) \in L^1(a,\infty)$
for a given $a>0$, for which there exists a constant $\a_f \ge 0$ such that
$$
f(x) \;\ll_\e\; x^{\a_f+\e}.\leqno(1.1)
$$
Actually it is more precise to define $\a_f$ as the infimum of the constants
for which (1.1) holds.
Here and later $\e > 0$ denotes arbitrarily small constants, not necessarily
the same ones at each occurrence. The notation  $A \ll_\e B$
(same as $A = O_\e(B)$) means that
$|A| \le C(\e)B$ for some positive constant $C(\e)$, which depends
only on $\e$. We define the convolution of functions $f,g\in {\Cal M}_1$ as
$$
(f\odot g)(x) := \int_1^x f(y)g({x\over y})\,{\d y\over y},\leqno(1.2)
$$
which is the special case $a=1$ of the more general convolution function
$$
(f\odot g)_a(x) := \int_a^{x/a} f(y)g({x\over y})\,{\d y\over y}
\qquad(a>0; f,g\in{\Cal M}_a).
\leqno(1.3)
$$
Of special interest is the function,
for $f\in\Cal M_1$,
$$
C[f(x)] := (f\odot f)(x)= \int_1^x f(y)f({x\over y})
{\d y\over y}\qquad(x\ge1),
\leqno(1.4)
$$
or more generally
$$
C_a[f(x)] := (f\odot f)_a(x)= \int_a^{x/a} f(y)f({x\over y})
{\d y\over y}\qquad(x\ge a, f\in {\Cal M}_a).
\leqno(1.5)
$$
Furthermore, the iterates of $C[f(x)]$ are defined as
$$
C^{(1)}[f(x)] \equiv C[f(x)],\quad C^{(k)}[f(x)] := C[C^{(k-1)}[f(x)]]
\quad(x\ge 1,
\, k \ge 2).\leqno(1.6)
$$
Obviously we have, in view of (1.1),
$$
C^{(k)}[f(x)] \ll_{\e,k} x^{\a_f+\e},\leqno(1.7)
$$
and in [10] the bound (1.7) was improved in case when $f(x)$
represents  several well-known number theoretic error terms.
In particular this includes the mean square and biquadrate of $|\zt|$
and the error terms in the corresponding asymptotic formulas,
$\D_k(x)$, the error term in the (generalized) Dirichlet
divisor problem and the problems involving the
distribution of non-isomorphic Abelian groups and the Rankin--Selberg
convolution of holomorphic cusp forms.
Relevant definitions and notions are to be found in [10].

\medskip
One of the reasons for the study of the convolution functions
$(f\odot g)(x)$ is that they appear naturally in the
context of (modified) Mellin transforms
$$
F^*(s) \equiv m[f(x)] := \int_1^\infty f(x)x^{-s}\d x
\qquad(s =\s + it;\; \s,t\in\RR)
\leqno(1.8)
$$
by means of the formula, which holds under suitable conditions,
$$
m[(f\odot g)(x)] \;=\; m[f(x)]\,m[g(x)].\leqno(1.9)
$$

\medskip
\medskip
The application to  the summatory function $A(x) := \sum_{n\le x}a_n$
of the sequence $\{a_n\}_{n=1}^\infty$ was given in [10].
Let $A(x)$ be of the form
$$
A(x) \;:=\; \sum_{i=1}^k\,\sum_{j=0}^{M_i}c_{i,j}x^{\a_i}\log^jx + u(x),
\leqno(1.10)
$$
where  the $c_{i,j}$'s are real constants with $c_{1,M_1} > 0$ and
$\a_1 > \a_2 > \ldots > \a_k > 0$,
and  $u(x)\;(= o(x^{\a_k})\,$ as $x\to\infty$)
is the error term in the asymptotic formula for $A(x)$.
If $u(x)$ satisfies  the mean square estimate
$$
\int_0^X u^2(x)\d x \;\ll\; X^{1+2\b}\qquad(0 \le \b < \a_k),
\leqno(1.11)
$$
then the following result was proved in [10].

\bigskip
THEOREM 1. {\it Let the above hypotheses on $A(x)$ and $u(x)$ hold,
and suppose that the function $\,{\Cal A}(s) = \sum_{n=1}^\infty
a_nn^{-s}\,$ admits analytic continuation
to the region $\R s > 0$, where it is regular except for the poles at
$s = \a_1, \a_2,\ldots, \a_k$ which are of order
$M_1+1, M_2+1,\ldots, M_k+1$, respectively. If}
$$
\int_T^{2T}|{\Cal A}(\s_1+it)|^2\d t \;\ll\; T^{2-\delta}
\leqno(1.12)
$$
{\it holds for some $\delta > 0$ and $0 < \s_1 < \a_k$, then we have}
$$
C[u(x)] \;=\; \int_1^x u(y)u({x\over y})\,{\d y\over y} \;\ll\; x^{\s_1}.
\leqno(1.13)
$$

\bigskip
\head
2. The asymptotics of the convolution function
\endhead

The asymptotic formula for $C[f(x)]$ is not easy to obtain,
even if a sharp formula for $f(x)$ (or its integral) is known. For example,
it is well known (see [5]) that
$$
\int_0^T|\zt|^2\d t = T\log\bigl({T\over 2\pi}\bigr) + (2\gamma-1)T
+ E(T),\leqno(2.1)
$$
where $\z(s) = \sum_{n=1}^\infty n^{-s}\;(\R s >1)$ is the Riemann
zeta-function, $\gamma = -\G'(1)$ is Euler's constant,
and for the error term $E(T)$ one has the asymptotic formula
$$
\int_0^TE^2(t)\d t = CT^{3/2} + O(T\log^4T)\qquad(C>0).\leqno(2.2)
$$
It seems difficult to obtain an asymptotic formula for $C[|\zx|^2]$,
even with the precise information contained in (2.1) and (2.2).
We shall return to this problem in Section 4.

In number theory one often encounters, as error terms
in asymptotic formulas, {\it regularly varying} functions. These
are  functions $h(x)$  which are positive, continuous (or, more generally,
measurable) for $x \ge x_0 \,(>0)$, for which
there exists $\rho \in\RR$ (called the index of $h(x)$) such that
$$
\lim_{x\rightarrow \infty} \frac{h(cx)}{h(x)}=c^{\rho},
\quad \text{for all }\; c>0. \leqno(2.3)
$$
We shall denote the set of all regularly varying functions
by $\Cal R$. We shall also denote by ${\Cal L}$ the set of
{\it slowly varying} (or {\it slowly oscillating})
functions, namely those functions $h(x)$ in $\Cal R$
for which the index $\rho=0$. It is easy to show that
if $h\in {\Cal R}$, then there exists $L\in {\Cal L}$
such that $h(x)=x^{\rho}L(x)$, with $\rho$ being the index of $h$.
\medskip

For a comprehensive account of regularly varying functions
the reader is referred to the
monographs of Bingham et al. [1] and E. Seneta [22]. By a fundamental
result of J. Karamata [16], who founded the theory of regular variation,
the limit in (2.3) is uniform for $0 < a \le c \le b < \infty$
and any $0 < a < b$. This is known as the {\it uniform convergence
theorem}. It is used to show that any slowly
varying function $L(x)$ (for $x \ge x_0 \;(>0))$ is necessarily of the form
$$
L(x) = A(x)\exp\left(\int_{x_0}^x\eta(t){\d t\over t}\right),\;
\lim_{x\to\infty}A(x) = A > 0,\;   \lim_{x\to\infty}\eta(x) = 0,\leqno(2.4)
$$
so that $x^{-\e} \ll L(x) \ll x^\e$ always holds.
If $h(x) \in {\Cal R}$ with index $\rho$, then
$$
C[h(x)] = \int_1^x h(u)h\bigl({x\over u}\bigr){\d u\over u}
= x^\rho\int_1^x L(u)L\bigl({x\over u}\bigr){\d u\over u}\,,
$$
where $L(x)$ is a slowly varying function.
Hence the problem of the asymptotic evaluation of $C[h(x)]$ is in this case
reduced to the evaluation of
$$
C[L(x)] = \int_1^x L(u)L\bigl({x\over u}\bigr){\d u\over u}
\qquad(L(x)\in {\Cal L}).\leqno(2.5)
$$

\smallskip
In some cases it is possible to evaluate the integral in (2.5) explicitly,
but in the general case it is not an easy task. For example, let $L(x) =
(\log x)^\a$ with $\a > -1$ a given constant. Then we have,
with the change of variable $t = \log u/\log x$,
$$
\eqalign{
C[(\log x)^\a] &= \int_1^x (\log u)^\a(\log x - \log u)^\a {\d u\over u}\cr
&= (\log x)^{2\a+1}\int_0^1 t^\a(1-t)^\a \d t\cr
&= {\G^2(\a+1)\over \G(2\a+2)}(\log x)^{2\a+1}.\cr}\leqno(2.6)
$$
We note that in (2.6) the resulting function is again slowly varying.
This is also true in general, when we consider $C_a[h(x)]$ (cf. (1.5))
for sufficiently large $a$ (if (2.4) holds, then $a=x_0$ may be taken).
The result is

\bigskip
THEOREM 2. {\it If $h(x)\in {\Cal R}$ with index $\rho$,
then for sufficiently large $a$ we have}
$$
C_a[h(x)] \;=\; x^\rho \int_a^{x/a} L(u)L\bigl({x\over u}\bigr){\d u\over u}
\;=\; x^\rho C_a[L(x)],\leqno(2.7)
$$
{\it where $C_a[L(x)]$ is a slowly varying function}.

\bigskip
{\bf Proof}.
Let $L(x) \in {\Cal L}, \,a>  0$. The result follows from
 the uniform convergence theorem that, uniformly for
$k_1 \le c \le k_2$, $L(x) \in {\Cal L}$ and any $k_2>k_1>0$ we have
$L(cx)/L(x)\to1$ as $x\to\infty$. Namely, if $B>0$ is a large constant, then
we have
$$
\eqalign{
C_a[L(x)] &= \int_a^{x/a} L(u)L\bigl({x\over u}\bigr){\d u\over u}
\ge \int_{\sqrt{x}/B}^{\sqrt{x}} L(u)L\bigl({x\over u}\bigr)
{\d u\over u}\cr&
= (1+o(1))L^2(\sqrt{x}\,)\int_{\sqrt{x}/B}^{\sqrt{x}}{\d u\over u}
\ge{\log B\over2}L^2(\sqrt{x}\,),\cr}\leqno(2.8)
$$
so that
$$
L^2(\sqrt{x}\,) \;\ll\;C_a[L(x)]/\log B\leqno(2.9)
$$
for $x\ge x(B)$. On the other hand, if $c\ge1$ is a given constant, then
$$
\eqalign{
C_a[L(cx)] &= 2\int_{a}^{\sqrt{cx}} L(u)L\bigl({cx\over u}\bigr){\d u\over u}
%\cr&
= (2+o(1))\int_{a}^{\sqrt{cx}} L(u)L\bigl({x\over u}\bigr){\d u\over u}\cr&
= (1+o(1))C_a[L(x)] + (2+o(1))\int_{\sqrt{x}}^{\sqrt{cx}}
L(u)L\bigl({x\over u}\bigr){\d u\over u}\cr&
= (1+o(1))C_a[L(x)] + O(L^2(\sqrt{x}))\int_{\sqrt{x}}^{\sqrt{cx}}{\d u\over u}\cr&
= (1+o(1))C_a[L(x)] + O(C_a[L(x)]/\log B),\cr}\leqno(2.10)
$$
since (2.9) holds. But $B$ can be arbitrarily large,
and consequently (2.10) implies that $C_a[L(cx)] \sim C_a[L(x)]$
as $x\to\infty$. This proves the assertion, since $c\ge1$ may be assumed without
loss of generality.

\head
3. Mean square bounds
\endhead
In case when it is difficult to obtain an asymptotic formula
for $C_a[f(x)]$ one has to be content with upper bound estimates.
In this direction we have

\bigskip
THEOREM 3. {\it Suppose that $f\in L^2(1,\,\infty)$ and that
for some $\;\t\ge0$ and $D \ge 0$ we have}
$$
\int_1^Xf^2(x)\d x \ll X^{1+2\t}(\log X)^D.\leqno(3.1)
$$
{\it Then, for any $a \ge 1$,}
$$
C_a[f(x)] \ll x^{\t}(\log x)^{c(\t)},\quad c(\t) = D+1 \;\;{\roman {if}}\;\;
\t> 0, \;\;c(\t) = D+2 \;\;{\roman {if}}\;\; \t=0.\leqno(3.2)
$$

\bigskip
{\bf Proof}.
We note that
$$
\eqalign{
C_a[f(x)] &= \int_a^{\sqrt{x}} f(u)f\bigl({x\over u}\bigr)\,{\d u\over u}
+ \int_{\sqrt{x}}^{x/a}f(u)f\bigl({x\over u}\bigr)\,{\d u\over u}\cr&
= 2\int_{\sqrt{x}}^{x/a}f(u)f\bigl({x\over u}\bigr)\,{\d u\over u}.\cr}
\leqno(3.3)
$$
The last integral is split into $\ll \log x$ subintegrals of the form
$$
I(x,T) := \int_T^{T'}f(u)f\bigl({x\over u}\bigr)\,{\d u\over u}
\qquad(\sqrt{x} \le T < T' \le 2T \le 2x/a).
$$
An application of the Cauchy-Schwarz inequality for integrals and (3.1)
gives
$$
\eqalign{
I(x,T) &\le \left(\int_T^{T'}f^2(u){\d u\over u}
\,\int_T^{T'}f^2\bigl({x\over u}\bigr){\d u\over u}\right)^{1/2}\cr&
= \left(\int_T^{T'}f^2(u){\d u\over u}\,\int_{x/T'}^{x/T}
f^2(u){\d u\over u}\right)^{1/2}\cr&
\ll (T^{2\t}(\log T)^{c(\t)-1}
\cdot (x/T)^{2\t}(\log x/T)^{c(\t)-1})^{1/2} \ll x^{\t}(\log x)^{c(\t)-1},
\cr}
$$
from which (3.2) follows.

\medskip
As an application of Theorem 3, we consider the distribution
of non-isomorphic Abelian groups. As usual, let $a(n)$ denote the number
of  non-isomorphic Abelian groups with $n$ elements (see e.g.,
[5, Section 14.5] for an extensive account). This is a multiplicative function
and its generating series is
$$
\sum_{n=1}^\infty a(n)n^{-s} = \z(s)\z(2s)\z(3s)\cdots\qquad(\R s >1).
$$
If one sets (this is (1.10) with $k=6,\;M_i \equiv 0$)
$$
A(x) := \sum_{n\le x}a(n) = \sum_{j=1}^6A_jx^{1/j} + R_0(x),
\quad A_j := \prod_{k=1,k\not=j}^\infty \z(k/j),
$$
then $R_0(x)$ can be thought of as the error term in the asymptotic
formula for the summatory function of $a(n)$. The author obtained
$$
C[R_0(x)] \;\ll_\e\; x^{1/6+\e}\leqno(3.3)
$$
in [10] directly, by using results on power moments of $\z(s)$.
A slight improvement of (3.3) follows from Theorem 3 (with $\t=1/6,\,D= 89$)
and the bound
$$
\int_1^X R_0^2(x)\d x \ll X^{4/3}(\log X)^{89},\leqno(3.4)
$$
of D.R. Heath-Brown [3], namely
$$
C[R_0(x)] \;\ll\; x^{1/6}(\log x)^{90}.
$$
Incidentally, the bound in (3.4) is best possible up to a
power of the logarithm, since the author [6] proved that
$$
\int_1^X R_0^2(x)\d x   \;=\;\Omega(X^{4/3}\log X),
$$
where as usual $f(x) = \Omega(g(x))$ means that $f(x) = o(g(x))$ does
not hold as $x\to\infty$. It seems reasonable to conjecture that
$$
C[R_0(x)] \;\sim\; x^{\rho}L(x) \qquad(L(x)\in{\Cal L},
\,0\le\rho\le 1/6,\;x\to\infty).
$$

\bigskip

\head
4. The second and fourth fourth power of the zeta-function
\endhead

\bigskip
In this section we shall consider the asymptotic evaluation of
$C[|\zx|^{2k}]$ when $k = 1$ and $k = 2$. Naturally, the values
$k > 2$ could be also considered, but the problem then becomes much
more difficult, since our knowledge on the $2k$-th moment of $|\zx|$
when $k>2$ is rather modest (see [5, Chapter 8]).

\medskip
It was proved in [10] that
$$
C[|\zx|^{2}] \;\ll_\e\; x^\e \leqno(4.1 )
$$
 and that
$$
C[|\zx|^4] \;\ll_\e\; x^\e\min(x^{4\mu({1\over2})},\,
x^{2\mu({1\over2})+{1\over4}},\, x^{1\over3}),\leqno(4.2)
$$
where for a given $\s\in\RR$
$$
\mu(\s) = \limsup_{t\to\infty} {\log|\z(\s+it)|\over\log t}\leqno(4.3)
$$
is the Lindel\"of function. If the famous (hitherto unproved)
Lindel\"of conjecture that $\mu(\s) = 0$ for $\s\ge\hf$ (or equivalently
that $\zt \ll_\e |t|^\e$) is true, then we have trivially
$$
C[|\zx|^{2k}] \;\ll_{\e,k}\;x^{\e}, \leqno(4.4)
$$
and in any case
$$
C[|\zx|^{2k}] \;\ll_{\e,k}\;x^{2k\mu({1\over2})+\e} \leqno(4.5)
$$
does hold. Heuristically, one expects $C[|\zx|^{2k}] \sim L(x)\;(\in{\Cal L})$
to hold as $x\to\infty$. More precisely, I conjecture that for $k = 1,2$
there exists a constant $A_k>0$ such that
$$
C[|\zx|^{2k}] \;\sim\; A_k(\log x)^{2k^2+1}\qquad(x\to\infty),\leqno(4.6)
$$
and (4.6) probably also holds at least for $k = 3$ and $k = 4$.
If true, this conjecture is certainly beyond reach at present. The heuristic
motivation for (4.6) is given shortly as follows.
For $k \ge 1$ a fixed integer let
$$
\int_0^T\vert\zeta(\txt{1\over 2} + it)\vert^{2k}\d t =
T\,P_{k^2}(\log T) + E_k(T),\leqno(4.7)
$$
where for some suitable coefficients $a_{j,k}\;(a_{k^2,k} >0)$ one has
$$
P_{k^2}(y) = \sum_{j=0}^{k^2}a_{j,k}y^j, \leqno(4.8)
$$
and in particular it is known  that $
P_1(y) = y + 2\gamma - 1 - \log(2\pi)$ holds (cf. (2.1)).
One hopes that
$$
E_k(T) = o(T) \qquad (T \to \infty)\leqno(4.9)
$$
will hold for every fixed integer $k \ge 1$, but so far this is
known to be true only in the cases $k = 1$ and $k = 2$, when
$E_k(T)$ is a true error term in (4.7) (see [5] and [8]).
Recently (see Conrey et al. [2]) plausible
heuristic arguments have been given, by employing the techniques of
random matrix theory, to produce explicit values of the coefficients
$a_{j,k}$ in (4.8). Nevertheless, the author in [8] expressed doubts
that (4.7)--(4.8) will, in general, hold for $k>4$. Regardless
of the moment conjecture, it certainly
seems plausible that, for some index $\rho = \rho(k) \ge 0$, one has
$$
C[|\zx|^{2k}] \;\sim\; x^\rho L(x) \in {\Cal R}\qquad(x\to\infty).
\leqno(4.10)
$$
If (4.7)--(4.9) holds, then for $\s>1$ and some constants $d_{j,k}$
 we have
$$
\eqalign{
{\Cal Z}_k(s) &:= \int_1^\infty |\zx|^{2k}x^{-s}\d x
= \int_1^\infty (x\,P_{k^2}(\log x) + E_k(x))'x^{-s}\d x\cr&
= \sum_{j=0}^{k^2+1}{d_{j,k}\over(s-1)^j} + s\int_1^\infty
E_k(x)x^{-s-1}\d x.\cr}\leqno(4.11)
$$
Thus we obtain analytic continuation of the Mellin transform ${\Cal Z}_k(s)$
to the region $\s\ge1$ (at least). From (1.9) it follows that
$$
C[|\zx|^{2k}] = {1\over2\pi i}\int_{1+\e-i\infty}^{1+\e+i\infty}
{\Cal Z}_k^2(s)x^{s-1}\d s.\leqno(4.12)
$$
We shift the line of integration in (4.12) to $\R s = c$ for some suitable
$0<c<1$, passing over the pole of ${\Cal Z}_k^2(s)$ of order $2k^2+2$.
By the residue theorem (4.6) follows, provided of course that we can
make this procedure rigorous.

\medskip
By the method of proof of Theorem 3 and (4.7)--(4.9) with $k = 2$ one
can easily improve (4.1) to
$$
C[|\zx|^2] \;\ll\; (\log x)^5.\leqno(4.13)
$$
Any further improvements seem difficult, but nevertheless we can
prove an asymptotic formula for the integral of $C[|\zx|^2]$, which
supports the conjectural (4.6) when $k=1$. This is

\bigskip
THEOREM 4. {\it There exist effectively computable constants $A \,(=1/6),
B,C,D$ such that}
$$
\int_1^XC[|\zx|^2]\d x = (A\log^3X + B\log^2X + C\log X + D)X + O_\e(X^{1/2+\e}).
\leqno(4.14)
$$

\bigskip {\bf Proof}. Integrating (4.12) when $k=1$ we obtain
$$
\int_1^XC[|\zx|^2]\d x = {1\over2\pi i}\int_{1+\e-i\infty}^{1+\e+i\infty}
{\Cal Z}_1^2(s){X^{s} - 1\over s}\,\d s.\leqno(4.15)
$$
We note (see the author's paper [11]) that the function $\Z_1(s)$
continues meromorphically to $\CC$, having only a double pole
at $s=1$, and  simple poles at $s = -1,-3,\ldots\;$. The principal
part of its Laurent expansion at $s=1$ is
$$
{1\over (s-1)^{2}} + {2\gamma - \log(2\pi)\over s - 1}.
$$
In (4.15) we shift the
line of integration to $\R s = \hf + \e$, passing over the pole $s=1$ of the
integrand of order four. By the residue theorem, the main term in (4.14)
comes from this pole. The integral over the line $\R s = \hf + \e$
is $\ll_\e x^{1/2+\e}$, if one uses the mean square bound
$$
\int_1^T|{\Cal Z}_1(\s+it)|^2\d t \;\ll_\e\; T^{2-2\s+\e}\qquad(\hf \le \s \le 1),
$$
proved in [13] by M. Jutila, Y. Motohashi and the author. The value $A= 1/6$ easily
follows by calculating the residue at $s=1$ of the integrand in (4.15).
\bigskip
The function $C[|\zx|^4]$ is more difficult to deal with than $C[|\zx|^2]$.
The results that we obtain in this case are contained in the following

\bigskip
THEOREM 5. {\it We have}
$$
C[|\zx|^4] \;\ll_\e\; \min(x^{2\mu({1\over2})+\e},\,
 x^{1\over4}(\log x)^{23/2}),\leqno(4.16)
$$
{\it and with suitable constants $A_j\;(j=0,\ldots,9)$ we have}
$$
\int_1^XC[|\zx|^4]\d x = X\sum_{j=0}^9 A_j\log^jX + O_\e(X^{5/6+\e}).\leqno(4.17)
$$

\bigskip {\bf Proof}. First note that
$$\eqalign{
C[|\zx|^4] &= \int_1^x|\zt|^4|\z\bigl(\hf + i{x\over t}\bigr)|^4\,{\d t\over t}\cr&
\ll_\e \int_1^x\left|\zt|^2|\z\bigl(\hf + i{x\over t}\bigr)\right|^2
t^{2\mu(1/2)+\e}(x/t)^{2\mu(1/2)+\e}\,{\d t\over t}\cr&
= x^{2\mu(1/2)+\e}\int_1^x|\zt|^2\left|\z\bigl(\hf +
i{x\over t}\bigr)\right|^2
\,{\d t\over t}\cr& = x^{2\mu(1/2)+\e}C[|\zx|^{2}] \ll_\e x^{2\mu(1/2)+\e}\cr}
$$
because (4.1) holds. This establishes the first bound in (4.16). The second one
follows from Theorem 3 (with $\t = 1/4,\,D=21/2$) and the bound
$$\eqalign{
\int_0^T|\zt|^8\d t & = \int_0^T|\zt|^2|\zt|^6\d t\cr&
\le \left(\int_0^T|\zt|^4\d t\int_0^T|\zt|^{12}\d t\right)^{1/2}\cr&
\ll (T\log^4T\cdot T^2\log^{17}T)^{1/2} = T^{3/2}\log^{21/2}T,
\cr}
$$
where the well-known bounds (see e.g., [5, Chapter 8])
$$
\int_0^T|\zt|^4\d t \ll T\log^4T,\quad \int_0^T|\zt|^{12}\d t \ll T^2\log^{12}T
$$
were used.
Note  that  the sharpest known result at present (see M.N. Huxley [4]) is
$\mu(1/2) \le 32/205 = 0.156097\ldots\;$, hence unconditionally we have the bound
$$
C[|\zx|^4] \ll x^{1/4}(\log x)^{23/2}.\leqno(4.18)
$$
The  proof of (4.17) is analogous to the proof of (4.14).
Note that we have, similarly to (4.15),
$$
\int_1^XC[|\zx|^4]\d x = {1\over2\pi i}\int_{1+\e-i\infty}^{1+\e+i\infty}
{\Cal Z}_2^2(s){X^{s} - 1\over s}\,\d s,\leqno(4.19)
$$
where ${\Cal Z}_2(s)$ is given by (4.11) with $k=2$.
This function is regular for $\s>\hf$, except for pole $s=1$ of order five
(see [13]).
Moreover we have the mean square bound (see the author's paper [12])
$$
\int_1^{T}|\Z_2(\s+it)|^2\d t \ll_\e T^{{15-12\s\over5}+\e}
\quad({\txt{5\over6}} \le \s \le {\txt{5\over4}}).\leqno(4.20)
$$
Thus (4.17) follows if we shift the line of integration in (4.19) to
$\R s = 5/6+\e$ and use (4.20); the main term in (4.17) comes from the residue
of the integrand at $s=1$. One can show that $A_9 = 1/(2520\pi^2)$ and evaluate
also explicitly the remaining constants $A_j\,(j = 0,\ldots\,,8)$.

\medskip
If the eighth moment bound holds for $|\zt|$ (cf. Theorem 3 with $\t=\e$), then
the right-hand side of (4.18) can be replaced by $x^\e$. Moreover, in this case
the exponent in (4.20) will be $4-4\s+\e$ for $\hf  <\s \le 1$, giving the
exponent $3/4+\e$ in the error term in (4.17).
%\medskip
\head
5. The Rankin--Selberg problem
\endhead

%\bigskip
This work will be concluded by analyzing estimates of convolution functions
in the classical Rankin-Selberg problem. In this section we shall
make a digression and consider
the problem itself by means of a complex integration technique, while
mean square bounds will be dealt with in the last section.
The  Rankin-Selberg problem consists
of the estimation of the error term function
$$
\D(x) = \sum_{n\le x}c_n - Cx,\leqno(5.1)
$$
where the notation is as follows (see e.g., R.A. Rankin's monograph [18]).
Let $\varphi(z)$ be a holomorphic cusp form of weight $\kappa$ with
respect to the full modular group $SL(2,\ZZ)$, and denote by
$a(n)$ the $n$-th Fourier coefficient of $\varphi(z)$.
We suppose that $\varphi(z)$ is a normalized eigenfunction for the Hecke
operators $T(n)$, that is,
$  a(1)=1  $ and $  T(n)\varphi=a(n)\varphi $ for every $n \in \NN$.
In (5.1) $C>0$ is a suitable constant (see e.g., [14]
for its explicit expression), and $c_n$ is the convolution function
defined by
$$
c_{n}=n^{1-\kappa}\sum_{m^2 \mid n}m^{2(\kappa-1)}
\left|a\Bigl({n\over m^2}\Bigr)\right|^2.
$$
The classical Rankin-Selberg  bound of 1939 is
$$
\D(x) = O(x^{3/5}),\leqno(5.2)
$$
hitherto unimproved. In their works, done independently,
R.A. Rankin [17] derives (5.2) from a general result of
E. Landau, while A. Selberg [20] states the result with no proof.
We shall  estimate now $\D(x)$ by the complex integration technique.
The key fact in this approach is that, for $s = \s+it$ with $\s > 1$,
one has the decomposition
$$
Z(s) := \sum_{n=1}^\infty c_n n^{-s} = \z(s)\sum_{n=1}^\infty b_n n^{-s}
= \z(s)B(s),\leqno(5.3)
$$
say, where $B(s)$ belongs to the Selberg class of Dirichlet series
of degree three, and $B(s)$ is
holomorphic for $\R s >0$. This follows from G. Shimura [23]
(see also A. Sankaranarayanan [19], who used (5.3) to obtain
mean square bounds for $Z(s)$).
The coefficients $b_n$ satisfy $b_n \ll_\e n^\e$ (see [19],
actually the coefficients  $b_n$ are bounded by a log-power in mean
square, but this is not needed here). For the definition and
properties of the Selberg class of $L$--functions the reader is referred to
A. Selberg [21] and the survey paper of Kaczorowski--Perelli [15].

\smallskip
On using classical Perron's formula (see e.g., the Appendix
of [5]) and the convexity bound $Z(s) \ll_\e
|t|^{2-2\s+\e}\;(0 \le \s \le 1,\; |t| \ge 1)$, it follows that
$$
\D(x) = {1\over2\pi i}\int_{{1\over2}-iT}^{{1\over2}+iT}{Z(s)\over s}x^s\d s
+ O_\e\left(x^\e\left(x^{1/2} + {x\over T}\right)\right)\quad(1\ll T\ll x).
\leqno(5.4)
$$
If we suppose that
$$
\int_X^{2X}|B(
\hf+it)|^2\d t \ll_\e X^{\t+\e}\qquad(\t\ge 1),\leqno(5.5)
$$
and use the elementary fact (see [5, Chapter 8] for
the results on the moments of $|\zt|\,$) that
$$
\int_X^{2X}|\zt|^2\d t \ll X\log X,\leqno(5.6)
$$
then from (5.3)--(5.6) and the Cauchy-Schwarz inequality for integrals
we obtain
$$
\D(x) \ll_\e x^\e(x^{1/2}T^{\t/2-1/2} + xT^{-1}) \ll_\e x^{{\t\over
\t+1}+\e}\leqno(5.7)
$$
with $T = x^{1/(\t+1)}$. Thus we have proved the following

\bigskip
THEOREM 6. {\it If} (5.5) {\it holds, then we have}
$$
\D(x) \;\ll_\e \;x^{{\t\over\t+1}+\e}.
\leqno(5.8)
$$

\medskip

As $B(s)$ belongs to the Selberg
class of degree three,
then $B(\hf+it)$ in (5.5) can be written as a sum of two Dirichlet polynomials
(e.g., by the reflection principle discussed in [5, Chapter 4]),
each of length $\ll X^{3/2}$. Thus by
the mean value theorem for Dirichlet polynomials (op. cit.)
we have $\t \le 3/2$, giving (with unimportant $\e$)
the Rankin-Selberg bound $\D(x) \ll_\e x^{3/5+\e}$.
Clearly improvement will come from better values of $\t$.
Note that the best possible
value of $\t$ in (5.5) is $\t = 1$, which follows
from general results on Dirichlet series
(see e.g., [5, Chapter 9]). It gives $1/2+\e$ as the
exponent in the Rankin-Selberg
problem, which is the limit of the method (the author's
conjectural exponent $3/8+\e$ (see [7]) is out of reach). To
attain this improvement one faces
essentially the same problem as in proving the sixth moment for $|\zt|$,
namely $\int_0^T|\zt|^6\d t \ll_\e T^{1+\e}$. In fact the present
problem is even more difficult, because the
properties of the coefficients $b_n$ are
even less known than the properties of the divisor coefficients
$$
d_3(n) \= \sum_{abc=n;a,b,c\in\NN}1,
$$
generated by $\z^3(s)$, which occur in the investigations relating to the
sixth moment of $|\zt|$.
If we knew the analogue of the strongest sixth moment bound
$$
\int_0^T|\zt|^6\d t \;\ll\;T^{5/4}\log^CT\qquad(C>0),
$$
namely $\t = 5/4$ in (5.5), then (5.7) would yield $\D(x) \ll_\e x^{5/9+\e}$,
 improving substantially (5.2).

%\bigskip
\head
6. Mean square and convolution in the Rankin--Selberg problem
\endhead
%\bigskip
In [14] the explicit formula for $\D(x)$ was derived. This is
$$
\D(x) = {x^{3/8}\over2\pi}\sum_{k\le K}c_kk^{-5/8}\sin
\left(8\pi(kx)^{1/4}+{3\pi\over4}\right) + O_\e(x^\e((Kx)^{1/4} + x^{3/4}K^{-1/4})),
\leqno(6.1)
$$
where $K$ is a parameter which satisfies $1 \ll K \ll x$.

If we use (6.1) with $K=x$, square and integrate, then by the first derivative test
(see e.g., [5, Lemma 2.1]) it follows that
$$
\int_1^X\D^2(x)\d x \;\ll_\e\; X^{1+2\b+\e}\leqno(6.2)
$$
holds with $\b = 1/2$. But as we have (see [14, eq. (3.5)])
$$
\D(X) = H^{-1}\int_{X-H}^{X+H}\D(x)\d x + O(H)\qquad(X^\e\le H \le \hf X),
\leqno(6.3)
$$
it follows by the Cauchy-Schwarz inequality that
$$
\D^2(X) \;\ll\; H^{-1}\int_{X-H}^{X+H}\D^2(x)\d x + H^2\qquad(X^\e\le H \le \hf X).
\leqno(6.4)
$$
Hence (6.2) with $\b = 1/2$ and (6.4) give (5.2) with the (poor)
exponent $2/3+\e$, and
any exponent $\b < 2/5$ would lead to an improvement of the Rankin-Selberg
exponent 3/5. Although we cannot at present attain such an improvement from
a mean square bound, we can improve on the value $\b = 1/2$. Namely, let as before
$\mu(\s)$ denote the Lindel\"of function (see (4.3)).
Then we have the following

\bigskip
THEOREM 7. {\it We have} (6.2) {\it with}
$$
\b \= {2\over 5-4\mu(\hf)}.\leqno(6.5)
$$

\bigskip
{\bf Proof}. From the analogy with the divisor problem (see e.g., [5, Chapter 13])
it follows that (6.5) will be proved if we can show that
$$
\int_T^{2T}|Z(\s+it)|^2\d t \;\ll\; T^{2-\delta}\leqno(6.6)
$$
holds with $\s > {2\over 5-4\mu(\hf)}$ and some small $\delta \,(>0)$,
with $Z(s)$ given by (5.3). Note that we have the functional equation
$$
Z(s) \= {\Cal X}(s)Z(1-s),\quad {\Cal X}(\s+it) \asymp |t|^{2-4\s}
\quad(0 < \s < 1),\leqno(6.7)
$$
since $Z(s)$ is in the Selberg class of degree four.
Furthermore, we have the mean square bound, proved by the author in [10, eq. (9.27)]
(in [10] the exponent of $T$ should have $4\mu(1/2)$ instead of $2\mu(1/2)$),
$$
\int_T^{2T}|Z(\s+it)|^2\d t \ll_\e\; T^{4\mu(1/2)(1-\s)+\e}(T+ T^{3(1-\s)})
\qquad(\hf\le \s \le 1).\leqno(6.8)
$$
Therefore we obtain
$$\eqalign{&
\int_T^{2T}|Z(\s+it)|^2\d t \ll T^{4-8\s}\int_T^{2T}|Z(1-\s+it)|^2\d t \cr&
\ll_\e T^{4-8\s+4\mu(1/2)\s+\e}(T + T^{3\s})\qquad(0< \s \le \hf),\cr}
$$
and for $1/3 \le \s \le 1/2$ the last quantity is $\ll T^{2-\delta}$
if $\s = (2+\delta+\e)/(5-4\mu(\hf))$, proving the assertion of Theorem 7.
Note that with the sharpest result (see M.N. Huxley [4])
$\mu(1/2) \le 32/205$ we obtain $\b = 410/897 = 0.4570709\ldots\,$.
The limit is the value $\b = 2/5$ if the Lindel\"of hypothesis
(that $\mu(\hf) =0$) is true.
Of course, improving the value $\t = 3/2$ in (5.5) would be another way
to improve on the value of $\b$.

\medskip
The merit of the value of $\b$ in (6.5) is that is strictly less than one half.
As already mentioned, if we square out and integrate (6.1),
all that follows is $\b \le\hf$. Incidentally, this bound follows
in the general case of the mean square bound for an $L$-function of degree four
in the Selberg class. Thus  Theorem 7 shows that the finer
information  that  we have
in the Rankin-Selberg problem (the product representation (5.3)) can be put
to advantage. As a consequence of Theorem 7 and Theorem 3 we obtain that
$$
C[\D(x)] \;\ll_\e\; x^{{2\over5-4\mu(1/2)}+\e}. \leqno(6.9)
$$
The bound (6.9) was obtained in [10] by a direct, more involved technique.
With some more effort one can replace `$\e$' in (6.9) by an explicit power
of the logarithm. If one considers averages of $C[\D(x)]$, then even more
cancellations occur. In this direction we shall prove

\medskip
THEOREM 8. {\it For any given $\e>0$ we have}
$$
\int_1^XC[\D(x)]\d x \;\ll_\e\; X^{5/4+\e}.\leqno(6.10)
$$
\medskip
{\bf Proof}. From (5.1) and (5.3) we obtain, for $\R s>1$,
$$
\eqalign{Z(s) &= \int_{1-0}^\infty x^{-s}\d\Bigl(\sum_{n\le x}c_n\Bigr)
= \int_{1-0}^\infty x^{-s}\Bigl(C\d x + \d\D(x)\Bigr)\cr&
= {Cs\over s-1} + s\int_1^\infty\D(x)x^{-s-1}\d x,\cr}
$$
since $\D(1-0) = -C$.
From (1.9) it follows that
$$
C[\D(x)] = {1\over2\pi i}\int_{1-i\infty}^{1+i\infty}{U^2(s)\over s^2}x^s\d s,
$$
where (5.3) shows that the function
$$
U(s) := Z(s) - {Cs\over s-1}
$$
is regular in the region  $\R s>0$. By integration we have
$$
\int_1^XC[\D(x)]\d x = {1\over2\pi i}\int_{1-i\infty}^{1+i\infty}{U^2(s)\over s^2}
\cdot {X^{s+1}-1\over s+1}\,\d s.\leqno(6.11)
$$
Now we shift the line of integration in the last integral to the line
$\R s = {1\over4}+ \e$.
We note that (6.7) holds, and we obtain that the right-hand side of (6.11) is
$$
\ll_\e X^{5/4+\e}\left(1 + \int_{-\infty}^{\infty}(|t|+1)^{-1-8\e}
|Z({\txt{3\over4}}-\e+it)|^2\d t\right) \ll_\e X^{5/4+\e}.\leqno(6.12)
$$
Namely $Z(s)$ is of degree four in the Selberg class,
and consequently by (6.7) and
the mean value theorem for Dirichlet polynomials one obtains without difficulty
$$
\int_1^T|Z(\s + it)|^2\d t \ll_\e T^\e(T + T^{4-4\s})\qquad(\hf \le \s \le 1).
\leqno(6.13)
$$
Then we obtain ($T \gg 1$), on using (6.13),
$$
\int_T^{2T}(|t|+1)^{-1-8\e}
|Z({\txt{3\over4}}-\e+it)|^2\d t \ll_\e T^{-1-8\e}T^{1+ 5\e} = T^{-3\e},
$$
which means that the integral in (6.12) converges, and (6.10) follows.
Finally we note that (6.13) can be sharpened to an asymptotic formula which
improves Theorem 3 of the author's paper [9]. This is

\bigskip
THEOREM 9. {\it If $\b$ is given by } (6.5), {\it then for fixed $\s$ satisfying
$\hf < \s \le 1$ we have}
$$
\int_1^{T}|Z(\s+it)|^2\d t = T\sum_{n=1}^\infty c_n^2n^{-2\s}
+ O_\e(T^{(2-2\s)/(1-\beta)+\e}).\leqno(6.14)
$$

\medskip
{\bf Proof}. We proceed as in the proof of Theorem 3 of [9]. The only difference
is that, instead of using (p. 174 of [9]) the bound
$$
\int_T^{2T}|E|^2\d t \,\ll_\e\, X^{2-2\s+\e} + T^2X^{1-2\s+\e},
$$
which corresponds to (6.2) with $\b = \hf$, we can use a better bound.
This is (6.2) with $\b$ given by (6.5), so that the above bound becomes
$$
\int_T^{2T}|E|^2\d t \,\ll_\e\, X^{2-2\s+\e} + T^2X^{2\b-2\s+\e},
$$
where $\b$ is given by (6.5) and satisfies ${2\over5} \le \b < \hf$. Instead
of the exponent $4-4\s+\e$ that appears in (4.2) of [9], we obtain now
the better exponent $(2-2\s)/(1-\b)+\e$ in (6.14). This ends the discussion
on Theorem 9, with the remark that its use instead of (6.13) does not lead
to a better exponent on the right-hand side of (6.10).

%\vfill
%\eject
%\topskip2cm

\bigskip\bigskip
\Refs
\bigskip\bigskip

\item{[1]} N.H. Bingham, C.M. Goldie and J.L. Teugels,
Regular Variation, CUP, Cambridge, 1987.

\item{[2]} J.B. Conrey, D.W. Farmer, J.P. Keating, M.O. Rubinstein
and N.C. Snaith, Integral moments of $L$--functions, Proc. London
Math. Soc. (3) {\bf91}(2005), 33-104.

\item{[3]} D.R. Heath-Brown, The number of Abelian groups of order
at most $x$, Ast\'erisque {\bf198-199-200}(1991), 153-163.

\item{[4]} M.N. Huxley, Exponential sums and the Riemann zeta-function V,
Proc. London Math. Soc. (3) {\bf90}(2005), 1-41.

\item{[5]} A. Ivi\'c, The Riemann zeta-function, John Wiley \&
Sons, New York, 1985 (2nd ed., Dover, Mineola, N.Y., 2003).

\item{[6]} { A. Ivi\'c}, The general divisor problem, J. Number Theory
{\bf 26}(1987), 73-91.

\item {[7]} { A. Ivi\'c}, Large values of certain number-theoretic
error terms, Acta Arith. {\bf56}(1990), 135-159.

\item{[8]} A. Ivi\'c, On some results concerning the Riemann Hypothesis, in
``Analytic Number Theory" (Kyoto,
1996) ed. Y. Motohashi, LMS LNS {\bf247}, Cambridge University Press,
Cambridge, 1997, pp. 139-167.

\item{[9]} A. Ivi\'c, On mean values of zeta-functions in the critical strip,
J. Th\'eorie des Nombres de Bordeaux {\bf15}(2003), 163-173.

\item{[10]} A. Ivi\'c,
Estimates of convolutions of certain number-theoretic error terms,
Inter. J. of Math. and Mathematical Sciences
2004:1, 1-23.

\item{[11]} A. Ivi\'c, The Mellin transform of the square of Riemann's
zeta-function, International J. of Number Theory  {\bf1}(2005),
65-73.

\item{[12]} A. Ivi\'c, On the estimation of some Mellin
transforms connected with the fourth moment of $|\zt|$, in Proc.
ELAZ2004 Conf. (Mainz, 2004), ed. W. Schwarz, in press,
{\tt ArXiv:math.NT \ /0404524}.

\item {[13]} A. Ivi\'c, M. Jutila and Y. Motohashi, The Mellin
transform of powers of  the Riemann zeta-function, {\it Acta Arith.}
{\bf95}(2000), 305-342.

\item{[14]} A. Ivi\'c, K. Matsumoto and Y. Tanigawa, On Riesz mean
of the coefficients of the Rankin--Selberg series, Math. Proc.
Camb. Phil. Soc. {\bf127}(1999), 117-131.

\item{[15]} A. Kaczorowski and A. Perelli, The Selberg class: a
survey, in ``Number Theory in Progress, Proc. Conf. in honour
of A. Schinzel (K. Gy\"ory et al. eds)", de Gruyter,
Berlin, 1999, pp. 953-992.

\item{[16]} J. Karamata, Sur un mode de croissance r\'eguli\`ere des
fonctions, Mathematica (Cluj) {\bf4}(1930), 38-53.

\item{[17]} R.A. Rankin,  Contributions to the theory of Ramanujan's
function $\tau(n)$ and similar arithmetical functions. II, The order
of Fourier coefficients of integral modular forms,  Math. Proc.
Cambridge  Phil. Soc. {\bf 35}(1939), 357-372.

\item{[18]} R.A. Rankin,  Modular Forms, Ellis Horwood Ltd.,
Chichester, England, 1984.

\item{[19]} A. Sankaranarayanan, Fundamental properties of symmetric
square $L$-functions I, Illinois J. Math. {\bf46}(2002), 23-43.

\item{[20]}A.~Selberg,  Bemerkungen \"{u}ber eine Dirichletsche Reihe,
   die mit der Theorie der Modulformen nahe verbunden ist,
 Arch. Math. Naturvid. {\bf 43}(1940), 47-50.

\item{[21]} A.~Selberg, Old and new conjectures and results about a class
of Dirichlet series, in ``Proc. Amalfi Conf. Analytic Number Theory 1989
(E. Bombieri et al. eds.)",
University of Salerno, Salerno, 1992, pp. 367--385.

\item{[22]} E. Seneta, Regularly varying functions, LNM {\bf508}, Springer
Verlag, Berlin--Heidelberg--New York, 1976.

\item{[23]} G. Shimura, On the holomorphy of certain Dirichlet series,
Proc. London Math. Soc. {\bf31}(1975), 79-98.

\vskip1cm
\endRefs

\enddocument

\bye